\documentclass[a4paper,11pt]{article}
\usepackage{amsmath,english,german,latexsym,amssymb}
\usepackage{rotating}
\newtheorem{thm}{Theorem}
\newtheorem{prop}[thm]{Proposition}
\newtheorem{lemma}{Lemma}
\newtheorem{cor}{Corollary}
\newtheorem{main}{Main Theorem}

\newcommand{\proof}[1][]{{\it Proof#1: }}

\newcommand{\qed}[1][1cm]{\hspace*{\fill} $\Box$ \vspace{#1}}

\renewcommand{\P}{{{\mathbf P}^1}}

\newcommand{\CC}{{\mathbf C}}

\newcommand{\cc}{{\mathbf {C}}}

\newcommand{\rr}{{\mathbf R}}

\newcommand{\del}{\partial}

\newcommand{\pairing}{\langle@,@,@,,@!@!@!@!@!@!\rangle}

\newcommand{\Orth}{{O}}

\newcommand{\Ospk}{\Orth^\prime_k}

\renewcommand{\a}{\alpha}

\newcommand{\cg}{\gamma}
\newcommand{\m}{\mu}

\newcommand{\D}{\Delta}
\newcommand{\e}{\varepsilon}

\renewcommand{\l}{\lambda}

\newcommand{\s}{\sigma}

\newcommand{\inj}{\hookrightarrow}

\newcommand{\sinfty}{\s_{\!\infty}}
\newcommand{\szero}{\s_0}

\newcommand{\Fk}{{\bf F}_{2k}}
\newcommand{\fk}{{\bf F}_{\!k}}
\newcommand{\fami}{{\cal F}}
\newcommand{\famit}{\tilde{\fami}}
\newcommand{\famimu}{{\cal F}_\m}

\newcommand{\milnor}{{\cal U}}

\newcommand{\xfami}{{\cal X}}

\newcommand{\labell}[1]{\label{#1}}

\newcommand{\G}{\Gamma}

\newcommand{\inv}{^{-1}}

\begin{document}

\thispagestyle{plain}

\begin{center}
\Large\sc
Monodromy groups of regular elliptic surfaces\\[1.5cm]
\end{center}

{\noindent
Michael L\"onne}\\


Monodromy in analytic families of smooth complex surfaces yields groups of
isotopy classes of orientation preserving diffeomorphisms for each family member
$X$.
For all deformation classes of minimal elliptic surfaces with $p_g>q=0$, we
determine the monodromy group of a representative $X$, i.e.\ the group of
isometries of the intersection lattice
$L_X:=H_2/\!\:\!$\raisebox{-0.7mm}{torsion}
generated by the monodromy action of all families containing $X$.
To this end we construct families such that any isometry is in the
group generated by their monodromies or does not respect the invariance of the
canonical class or the spinor norm.


\section*
{Introduction}

Monodromy is a well exploited concept in the realms of complex curves and
complex singularities, \cite{sing}.
In connection with complex surfaces monodromy of Lefschetz fibrations was
employed for the investigation of elliptic surfaces in the eighties.
The results are concisely presented in \cite{fm}.
Quite recently interest was revived and new results on higher genus Lefschetz
fibrations have emerged, e.g.\ \cite{mats}, \cite{oz}.\\
After the surge of Donaldson theory it was noted that monodromy techniques applied
to complex surfaces can serve to show the smooth invariance of the divisibility of
the canonical class. This was used to present lots of examples of
non-diffeomorphic complex surfaces of equal topological type,
\cite{fmm},\cite{mois},\cite{sal}.\\
Especially for complete intersection surfaces with positive geometric genus it was
shown that the monodromy group is a group of index two only in the
group of orthogonal transformations fixing the canonical class \cite{habil} and of
index two also in the image of the
natural representation of the group of diffeomorphisms, \cite{eo2}.\\ 
It is the intention of this paper to extend this result -- which
shows a astonishing likeness between differential topology notions and analytical
ones as in the curve case -- to the important class of regular minimal elliptic
surfaces with positive geometric genus.\\
The paper is organized in three paragraphs. First for each deformation type of
elliptic surfaces with $p_g>q=0$ and multiplicities $\mu=m_1,...,m_n$ of multiple
fibres we construct a family $\famimu$. The associated monodromy is shown to be
given by a complete vanishing lattice, a sublattice of corank two of the
intersection lattice together with its classes of square $-2$, as the group
generated by the reflections corresponding to these classes.\\
In order to find more monodromy we then construct families which will allow to
include the classes of multiple fibres as additional generators of the
corresponding extended vanishing lattice.\\
Finally we employ results on the deformation equivalence of elliptic surfaces to
fit the monodromies together for a proof of the main theorem. 

\begin{main}
\labell{main}
Let $X$ be a minimal elliptic surface with positive geometric genus $p_g$ and
vanishing irregularity $q$, then there exist families of elliptic surfaces
containing $X$, such that the induced monodromy actions on the homology lattice
$L_X$ generate $O'_k(L_X)$, the group of isometries of real spinor norm one
fixing the canonical class.
\end{main}

\section*
{Families of elliptic surfaces as double covers}
\newcommand{\cctcc}{\cc\times\cc}

The basic idea is, that families of regular elliptic surfaces without multiple
fibres which contain many degenerations can be obtained as double covers of
trivial families of Hirzebruch surfaces provided a suitable family of branching
curves containing many singular ones has been chosen.\\
Any Hirzebruch surface $\fk,\,k\!>\!0$ is ruled over $\P$ with a unique section
$\sinfty$ such that the complement is isomorphic to the line bundle ${\cal
O}_\P(k)$, which has a zero section~$\s_0$.
Having assigned coordinates $(x_0\!:\!x)$ on $\P$, let $\ell$ denote the line
of the ruling mapping to $(0\!:\!1)$ and let the complement of $\sinfty\cup\ell$
be trivialized as $\cc\times\cc$ with fibre coordinate $y$ and base
coordinate $x$.\\
Given the function $g_0(x,y)=y^3-x^{3k-1}$ there is a linear function $l(x,y)$
such that the function $g_\e:=g_0+\e l$ is a Morsification of $g_0$ for almost
all sufficiently small $\e$.
The $u$-level set of a function $g_\e$ is then a curve in $\cctcc$.
Compactification of these curves inside $\fk$ yields a two parameter set
$\{D'_{u,\e}\}$ of compact curves. Furthermore let $D_{u,\e}$ denote the
curve decomposing into $D'_{u,\e}$ and the section $\sinfty$.

\begin{lemma}
\labell{lin-sys}
The curves $D_{u,\e}$ are linearly equivalent to $\sinfty\!+\!3\szero$ and
smooth in a regular
${\cal C}^\infty$-neighbourhood of $\ell\cup\sinfty$.
\end{lemma}

\proof
Consider $\cc\times\cc$ with coordinates $w,x_0$ and the curves defined by
\begin{equation}
\label{curve-eqn}
w^3-x_0+\e l(x_0^{2k}w,x_0^{3k-1})-ux_0^{3k}=0.
\end{equation}
Then this chart and its curves glue to the chart and its curves given above via
the identification
$$
w=x^{-k}y,x=x_0^{-1},
$$
and yield compactifications of the $u$-level sets inside ${\cal O}_\P(k)$ for
each pair of parameters.
Since they are covers of the base of degree $3$, curves in
$\{D'_{u,\e}\}$ are linearly equivalent to $3\s_0$ and the claim for the
equivalence class of a curve $D_{u,\e}$ is immediate.
Moreover curves in $\{D'_{u,\e}\}$ are disjoint from $\sinfty$ and intersect
$\ell$ in the point $(w,x_0)=(0,0)$. Hence a curve $D_{u,\e}$ is smooth
along the smooth section $\sinfty$ and along $x_0=0$ as can be read off equation
(\ref{curve-eqn}). Smoothness finally extends to some regular neighbourhood of
$\sinfty\cup\ell$.
\qed

Next we fix some open ball $B\subset\cc^2$ containing the origin.
The double cover $\famit$ of $\Fk\times B$ branched along $\{(p,u,\e)|p\in
D_{u,\e}\}$ then fibres naturally over $B$ with any smooth fibre $X$ a double
cover of $\Fk$ branched along some smooth divisor $D_{u,\e}$.
The ruling of $\Fk$, its distinguished fibre $\ell$ and section $\sinfty$ pull
back to an elliptic fibration, a cusp fibre and a section $\s$.
Computation of the numerical invariants yields
$$
p_g=k-1,q=0,\chi=k,{c_1}^2=0.
$$
and therefore the surfaces of the family are relatively minimal elliptic surfaces
with no multiple fibres.\\

Over the subset $\cc^2\times B$ the double cover has an explicit description as
the hypersurface in $\cc^3\times B$ given by the equation
$$
z^2+g_\e(x,y)-u=0.
$$
Let $\milnor$ be the open subset of $\famit$ given by the subset of this
hypersurface cut out by a big transversal ball in $\cc^5$.\\
Locally at the origin in $B$ we perform simultaneous logarithmic transforms, cf.\
\cite[p.112]{fm} along smooth fibres in the vicinity of the cusp fibre off
$\milnor$. So for each unordered finite sequence $\mu=m_1,...,m_n$ of
multiplicities we get a family
$\tilde{\famimu}$ over a possibly smaller $B_\mu$ containing the origin such
that a smooth fibre $X_\mu$ is an elliptic surface with invariants as above but
with $n$ multiple fibres of the given multiplicities.\\
If now $\e=\e_\mu$ is chosen sufficiently small and general, $g_\e$ is a Morse
function with all its critical values in $B_\mu$.
Hence the family $\famimu$ obtained by restriction of $\tilde\famimu$ to the
$\e$-slice of $B_\mu$ has the property, that $\famimu\cap\milnor=\fami\cap\milnor$
is a Morsification of the hypersurface singularity given by $g_0$ with Milnor fibre
$M=X\cap\milnor$.\\

Our final aim in this paragraph is to get hold on
the algebraic monodromy group of the family $\famimu$ and a set $\D$ of classes of
square $-2$ by which the monodromy is determined as the group $\G_\D$ generated
by the reflections on hyperplanes normal to elements in $\D$.
\begin{prop}
Let $X_\mu$ be a smooth surface of geometric genus $p_g\!>\!0$ in a
family~$\famimu$. Then the monodromy coincides with the group $\G_{\D_M}$
associated to the set of all classes of square $-2$ contained in the even
unimodular sublattice $L_M$ of corank two of classes supported on the Milnor fibre
$M$. Moreover the pair
$L_M,\D_M$ is a complete vanishing lattice in the sense that, cf.\
\cite[5.3.1]{habil}:
\begin{enumerate}
\item
$\D_M$ generates $L_M$,
\item
$\D_M$ is a single $\G_{\D_M}$ orbit,
\item
$\D_M$ contains six elements the intersection diagram of which is\\[6mm]
\unitlength1.6cm
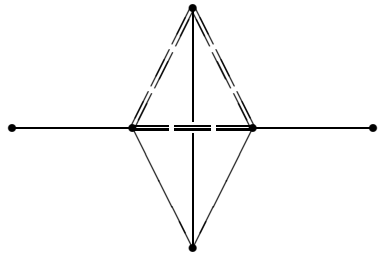
\begin{picture}(4,1.6)(-1,1.2)
\put(1,2){\circle*{.06}}
\put(1,2){\line(1,0){1}}
\put(2,2.015){\line(1,0){.3}}
\put(2.35,1.985){\line(1,0){.3}}
\put(2.35,2.015){\line(1,0){.3}}
\put(2.7,1.985){\line(1,0){.3}}
\put(2.7,2.015){\line(1,0){.3}}
\put(2,1.985){\line(1,0){.3}}
\put(3,2){\line(1,0){1}}
\put(2,2){\circle*{.06}}
\put(2.01,1.98){\line(1,2){.15}}
\put(1.98,2){\line(1,2){.15}}
\put(2.185,2.33){\line(1,2){.15}}
\put(2.155,2.35){\line(1,2){.15}}
\put(2.36,2.68){\line(1,2){.15}}
\put(2.33,2.7){\line(1,2){.15}}
\put(2,2){\line(1,-2){.5}}
\put(3,2){\circle*{.06}}
\put(2.845,2.35){\line(-1,2){.15}}
\put(2.815,2.33){\line(-1,2){.15}}
\put(2.67,2.7){\line(-1,2){.15}}
\put(2.64,2.68){\line(-1,2){.15}}
\put(3.02,2){\line(-1,2){.15}}
\put(2.99,1.98){\line(-1,2){.15}}
\put(3,2){\line(-1,-2){.5}}
\put(4,2){\circle*{.06}}
\put(2.5,3){\circle*{.06}}
\put(2.5,3){\line(0,-1){.95}}
\put(2.5,1){\circle*{.06}}
\put(2.5,1){\line(0,1){.95}}
\end{picture}
\end{enumerate}
\end{prop}

\proof
Since $\famimu$ is a family of surfaces with at most a single Morse
singularity, the monodromy group is generated by the reflections associated to the
vanishing cycles.
They coincide with the vanishing cycles of $\famimu\cap\milnor$ for outside
$\milnor$ the family projection is submersive by construction.
We made sure that $\famimu\cap\milnor$ is the Morsification of the hypersurface
singularity $z^2+y^3+x^{6k-1}$, $k\geq2$ of type ${\bf E}_{12k-4}$ in the list of
Arnold.\\
Hence the Milnor lattice, being even unimodular according to \cite[tab.3]{e-quad},
maps onto a sublattice $L_M$ of corank two in $L_{X_\mu}$ and is together with the
set of its classes of square $-2$ a complete vanishing lattice according to
\cite[proof of 4.1.2]{habil}.
\qed

\section*{Families with vanishing cycles around multiple fibres}

\newcommand{\dmu}{\D^\mu}

Given $p_g$ and $\mu=m_1,...,m_n$ we have considered a family $\famimu$ of regular
elliptic surfaces with geometric genus $p_g$ and multiple fibres $f_i$ of
multiplicities $m_i$.
A suitable restriction yields a family $\famimu|$ over a disc $\dmu$ with a single
quadratic surface degeneration at its origin.\\
We first shift our attention to the associated degeneration family of divisors on
$\Fk$.
In a suitable open $U\subset\Fk$ the singular divisor consists of two discs
which cross normally in a single point and map bijectively to some disc $\D$ in
the base of the ruling~$\pi:\fk\to\P$.
Moreover we may assume for $\l\in\dmu$ small that the intersection with $U$ of
the corresponding divisor $D_\l$ is a connected double cover of $\D$ with two branch
points except for $D_0$ branching at the origin only.\\
To each parameter $\l\in\dmu$ we then associate
$$
B_\l:=\{x\in\D|x\text{ is branch point of }D_{\l'}\to\D\text{ for some $\l'$ with
}|\l'|=|\l|\}
$$
which is a smooth one-cycle, degenerate in case $\l=0$ and enclosing the origin
for~$\l$ in a punctured neighbourhood of the origin.
The cycle $B_\l$ decomposes naturally into two arcs $B^{1,2}_\l$ connecting the
two branch points corresponding to $\l$.
To each radial ray in $\dmu$ we thus obtain two smooth families of arcs connecting
the branch points and degenerating to the origin of $\D$.
The complement in $\Fk$ of the section at infinity is a line bundle over $\P$ which
over $\D$ can be trivialized as $\D\times\CC$ of which~$U$ becomes a subset.
Having thus made sense of the straight segment $\overline{p_1p_2}$ between points
$p_1,p_2$ in $U$ we may consider
$$
\D^{1,2}_\l:=\{p\in\Fk|p\in\overline{p_1p_2},p_1,p_2\in D_\l\cap
U,\pi(p_1)=\pi(p_2)\in B^{1,2}_\l\}
$$
By pull back along the double cover $d_\l:X_\l\to\Fk$ branched at $D_\l$ we get
$$
S^{1,2}_\l:=\{p\in X_\l|d_\l(p)\in\D^{1,2}_\l\}.
$$

\begin{lemma}
The sets $S^{1,2}_\l$ are vanishing cycles for the surface degeneration over
$\dmu$.
\end{lemma}

\proof
Each arc $B^i_\l, i=1,2$ is lifted to a non-intersecting pair of arcs in $\Fk$ with
identical endpoints. The cycle thus obtained is filled to yield the disc
$\D_\l^{i}$ with boundary on the branch divisor $D_\l$.
This disc in turn is lifted to a pair of non-intersecting discs in $X_\l$ with
common boundary yielding the sphere $S_\l^i$, which contracts to a point when $\l$
tends to zero and the arc contracts to the origin of $\dmu$.
\qed

Next we consider holomorphic sections $\s^{1,2}$ to the projection
$\dmu\times\D\to\dmu$ which are affine linear in the second component
$$
\s_2^{1,2}(\l)=\pm x_0(\l_0-\l),\quad\text{ with }
x_0\in\D,\l_0\in\rr^{\geq0}\cap\dmu.
$$
For $\l$ varying along the non-negative reals the branch points move along two
real curves in $\D$ emerging from the origin with opposite tangent vectors.
For $\dmu$ suitably small there is thus a pair $\s^{1,2}$ as above such that
$\s_2^{i}(\l)$ is disjoint from $B^{j}_\l$ for $i\neq j$ and all $\l$ in an open
neighbourhood $W$ of the positive reals in $\dmu$.
Moreover there is a ray $\rho$ in $\D$ which intersects $B^1_{\l_0}$ transversally
once and is disjoint from $B^2_{\l_0}$.\\
In two copies of the family of elliptic surfaces over $W$ we cut off the
regular fibres which map to the image of $\s^1$, resp. $\s^2$. But still
$S^2_{\l_0}$, resp. $S^1_{\l_0}$ is a vanishing cycle by construction.
Moreover in the surface $X^\circ_\mu$ at the parameter $\l_0$ the same regular fibre
is cut off, hence we may glue both copies along this surface and thus get
a family $\xfami_\mu^\circ$ of elliptic surfaces degenerating at two points of its
base. The preceding family of the two copies we call $\xfami_\mu$.

\begin{lemma}
There are homology classes $\a,\a'$ for the surface $X^\circ_\mu$ which are
represented by vanishing cycles for the family $\xfami_\mu^\circ$ such that
$\a+\a'$ is a primitive class in the collar of $X^\circ_\mu$ orthogonal to a local
section.
\end{lemma}

\proof
Let $V:=f\inv(\D^*)$ be the preimage of $\D^*$ under the fibration map.
The collar $C_X$ of $X$ is homotopically equivalent to a real $3$-torus and can thus
be given as a handlebody with a single $3$-handle, three $2$-handles and three
$1$-handles. It embeds as a collar into $V$ which is given as a handlebody just by
adding two further $2$-handles with framing $-1$ along circles of the same isotopy
class.\\
We thus get that $H_3(V,C_X)=0$ and $H_2(V,C_X)$ is free abelian, and conclude
that $H_2(C_X)$ injects into $H_2(V)$ with torsionfree cokernel of rank one.
Since $H_2(V)$ contains elements of non-zero self-intersection but $H_2(C_X)$ only
elements of zero self-intersection, this injection is onto the radical of
$H_2(V)$.\\
Of course we aim at giving the cycles $S^{1,2}_{\l_0}$ in $X_\mu^{o}$ the right
orientations. In the corres\-ponding compact surface in $\xfami_\mu$ they
were vanishing cycles for the same
degeneration and thus could be given orientations such that the sum of their
classes became trivial. With the same orientations they represent classes $\a,\a'$
of $X_\mu^{o}$ the sum of which is still trivial though only numerically. Therefore
$\a\!+\!\a'$ is represented by a cycle in~$C_X$.\\
Finally notice that the fibration over the ray $\rho$ is trivial, and that
$S^{1,2}$ intersect this fibration in an embedded nontrivial circle on a single
fibre. Choosing a dual circle we get by moving
this circle along $\rho$ a relative class in $H_2(V,\del)$ which is dual to
$\a+\a'$. Having shown the sum to be primitive, it is
representable in the collar by a primitive element too, and by construction both
$\a,\a'$ are orthogonal to the local section which maps to the section at infinity
of the ruled surface.
\qed

In a sense to be made precise below instead of cutting off regular fibres we may
replace them by means of simultaneous logarithmic transformations.

\begin{prop}
\labell{xmfami}
For any $m$ there is an elliptic surface $X$ and classes $\a_m,\a'_m$ which are
represented by vanishing cycles with respect to some family $\xfami_{\mu,m}$
such that the class of a fibre of multiplicity $m$ is contained in a sublattice of
$L_X$ generated by $f,\a_m,\a'_m$.
\end{prop}

\proof
The family $\xfami$ and the section $\s^{1,2}$ constitute a situation as in
\cite[I.7.2]{fm}:
$$
\begin{array}{ccccccc}
\xfami & \to & \P\times W\vee W\\
\downarrow && \downarrow & \uparrow & \s^1\vee\s^2\\
W\vee W & = & W\vee W
\end{array}
$$
We get simultaneous logarithmic transformations locally over open subsets covering
a neighbourhood of the ray of positive reals,
having specified a choice of torsion point and a lift thereof.\\
The resulting family $\xfami$ has of course the degenerations which provide the
classes $\a,\a'$ of the previous lemma in the complement of the multiple fibre.\\
We infer from the discussion in \cite[2.1.2]{fm} of the topology of the logarithmic
transform that the second homology of its collar surjects onto the second homology
of a neighbourhood of a multiple fibre.
Hence the class of the multiple fibre is a linear combination of the fibre class
which has a single transversal intersection with a local section of the collar and
of a primitive class orthogonal to this section.
By the previous lemma a certain class of this kind is
the sum of the vanishing cycles. Since we are free to make our choice of torsion
point accordingly, we may assume these classes to coincide and get our claim.
\qed

\section*
{Monodromy groups}

Two minimal elliptic surfaces with positive Euler number are deformation
equivalent through elliptic surfaces if and only if their geometric genera,
their irregularities and their multiplicities of multiple fibres coincide.

This assertion is a theorem in \cite[II.7.6]{fm} given that the said
invariants determine the Euler number of a fibration and the diffeomorphism type
of the orbifold structure on the base and vice versa.

The notion of deformation equivalence through elliptic surfaces applies as soon
as an equivalence can be established via families with a compatible global
elliptic fibration structure -- as e.g.\ in our families $\famimu$.
It has the following important implication:

\begin{prop}
\labell{defo-fix}
Any orthogonal transformation of homology lattices induced from a deformation
equivalence through elliptic surfaces preserves the fibre class.
\end{prop}

Similarly there is the notion of deformation equivalence through elliptic surfaces
with a section which occurs if a compatible global section to the global elliptic
fibration can be given, here $\fami$ may serve as an example.
Again there is an associated invariance property:

\begin{prop}
\labell{defo-fixx}
Any orthogonal transformation of homology lattices induced from a deformation
equivalence through elliptic surfaces with a section preserves the fibre class and
the section class.
\end{prop}

In a situation that two given surfaces $X,X'$ are deformation equivalent and a
subgroup of the monodromy group of $X$ is known, it is an obvious goal to obtain
results for the monodromy group of $X'$.\\
To achieve this goal we take the family providing the monodromy and a family
providing the deformation equivalence and glue them by identification of the
surface $X$ only, which occurs in both, and its base points. Thus we may pull back
the monodromy group along an lattice isomorphism provided by the deformation
equivalence.

Our first application of this strategy exploits the fact \cite[p.57]{fm}, that
two elliptic surfaces with a section of equal irregularity and geometric genus are
deformation equivalent through elliptic surfaces with a section.

\begin{lemma}
\label{defo-sec}
For any given integer $p_g>0$ there is a minimal regular elliptic surface~$X$ of
geometric genus $p_g$ without multiple fibres and classes $\a_1,\a_1'$ in its
homology lattice, which are represented by vanishing cycles in suitable
degeneration families of~$X$, such that $\a_1+\a_1'$ is the fibre class.
\end{lemma}

\proof
Consider the divisor of type $4\s_0$ in $\fk$ given by the equation
$$
y^4-2y^2+\frac12=x^{4k}.
$$
The corresponding double cover is a smooth minimal regular elliptic surface with
$p_g=k-1$ which has fibres of type $I_2$ such that the pair of components of such
a fibre is dual to the pair of sections contained in the preimage of the section
at infinity of $\fk$.\\
After a choice of a section there is a deformation equivalence of $X$ through
elliptic surfaces with section to a smooth surface in $\fami$ by which the
component disjoint from the chosen section becomes a vanishing cycle for the
composed family. Thus the classes of both components can be represented by
vanishing cycles of appropriate families. Since their sum is the fibre class, we
are done.
\qed

We can make use of this result also in the presence of multiple fibres due again to
the existence of local simultaneous transformations.

\begin{prop}
\labell{fib-gen}
Given any integer $p_g>0$ and multiplicities $m_1,...,m_n$ there is an
elliptic surface $X_\mu$ with geometric genus $p_g$ and multiple fibres $f_i$ of
multiplicities $m_i$ and classes $\a_1,\a_1'$ in its homology lattice  which are
represented by vanishing cycles of appropriate degeneration families and the sum
of which is the class of a general fibre.
\end{prop}

\proof
By the preceding lemma there is a suitable surface in the absence of multiple
fibres. On the degeneration families providing the vanishing of the classes of the
components of an $I_2$ fibre we perform simultaneous logarithmic transforms to get
to our claim.
\qed

Another extension of the lemma is needed in the exceptional case of $K3$ surfaces.

\begin{prop}
\label{K3ex}
There is an elliptic $K3$ surface $X$ with classes $\a,\a',\s$ in its homology
lattice which are represented by vanishing cycles in appropriate degeneration
families such that
$$
\a+\a'=f, \a.\s=1,\a'.\s=0.
$$
\end{prop}

\proof
Let $X$ be the double cover of the quadric branched along the smooth branch curve
of bidegree $(4,4)$ with equation
$$
(x_1^2+x_0^2)^2(y_1^2+y_0^2)^2+(x_1^2+x_0^2)^2y_1^2y_0^2
+x_1^2x_0^2(y_1^2+y_0^2)^2-x_1^2x_0^2y_1^2y_0^2=0.
$$
Consider first the ruling of $\P\times\P$ which contains the line $x=\rho$
($\rho^3=1,\rho\neq1$), a bitangent to the branch curve.
The pull back along the double cover gives an elliptic fibration of $X$ with a
singular fibre of type $I_2$ at $x=\rho$. As in the lemma we get degeneration
families with vanishing cycles representing the homology classes $\a,\a'$ of the
fibre components.\\
The analogous result holds for the second ruling and the line $y=\rho$. The two
components of the corresponding fibre of type $I_2$ are dual to $\a,\a'$.
Therefore we may finish the proof by choosing $\s$ accordingly.
\qed 

We may now apply our pull back strategy to arbitrary regular elliptic surfaces:

\begin{prop}
\labell{v-pb}
Let $X$ be a minimal elliptic surface with $p_g>q=0$. Then the set $\D$ of classes
in $L_X$ of square $-2$ orthogonal to the fibre class $f$ and represented by a
vanishing cycle in a degeneration family of $X$ is such that
\begin{enumerate}
\item
there are classes $\a_1,\a_2\in \D$ with $\a_1+\a_1'=f$,
\item
for any class $f_m$ represented by a fibre of multiplicity $m$ there are classes
$\a_m,\a_m'\in\D$ with $f_m$ in the sublattice of $L_X$ generated by
$\a_m,\a_m',f$,
\item
there is a complete vanishing lattice $L',\D'$ with $\D'$ contained in $\D$ and
$L'$ an even unimodular sublattice of $L_X$ of corank two.
\end{enumerate}
\end{prop}

\proof
By the numerical criterion stated above $X$ is deformation equivalent through
elliptic surfaces to some surface $X_\mu$ as given in prop.\,\ref{fib-gen} and
further to a smooth surface in $\famimu$.
The vanishing cycles for these pull back to
elements of $\D$ such that the properties as in the first and third claim are
established.\\
If $\mu=m_1,...,m_n$ are the multiplicities of $X$ then $X$ is deformation
equivalent through elliptic surfaces to a smooth surface in an appropriate family
$\xfami_{\hat\mu,m}$, with $\hat\mu=m_1,...\hat m,...,m_n$, i.e.\ $m$ omitted once.
The pull back of the associated vanishing classes as given in prop.\,\ref{xmfami}
provide $\a_m,\a_m'$ such that a class of a fibre of multi\-plicity~$m$ is
in the sublattice generated by $\a,\a',f$.
\qed 

\begin{cor}
\labell{D-span}
Given $X$ and $\D$ as above, the orthogonal complement $L$ of the fibre class in
$L_X$ is generated by elements of $\D$.
\end{cor}

\proof
The lattice $L'$ of the proposition is primitive in $L$ of corank one. The radical
of $L$ is non-trivial and hence of rank one and a direct sum complement to $L'$.
It is generated by the classes of a general and of multiple fibres.\\
We may thus argue as follows. $L'$ is
generated by its $-2$ classes all of which belong to $\D$ since $L'$ is a
vanishing lattice.
The fibre class is contained in the sublattice of $L$ generated by
$\D$ by the second assertion of the proposition.
Finally also the classes of multiple fibres belong to the linear span of the
elements of
$\D$ by the second assertion.
\qed

\begin{cor}
\labell{D-conj}
Given $X$ and $\D$ as above, the elements of $\D$ are conjugate under the action
of $\G_\D$.
\end{cor}

\proof
Since $L',\D'$ is a complete vanishing lattice the elements of $\D'$ are
contained in a single orbit of the $\G_\D$ action.\\
Furthermore any $\a\in\D$ differs only by an element of the radical of $L$ from an
element $\a^*\in L$. Thus $\a^*$ is in $\D'$ and we argue again with the fact that
$L',\D'$ is a complete vanishing lattice to get an element $\delta\in\D$ with
$\a^*.\delta=1=\a.\delta$.
It is easy to check that $\cg_\delta\circ\cg_\a(\delta)=\a$.
\qed

\begin{cor}
\labell{D-special}
Given $X$ and $\D$ as above, there are six elements of $\D$
with intersection pattern given by\\[2mm]
\unitlength1.6cm
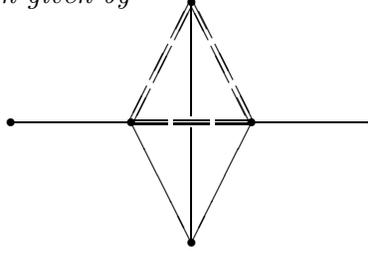
\begin{picture}(4,2.8)(-1,0)
\put(1,2){\circle*{.06}}
\put(1,2){\line(1,0){1}}
\put(2,2.015){\line(1,0){.3}}
\put(2.35,1.985){\line(1,0){.3}}
\put(2.35,2.015){\line(1,0){.3}}
\put(2.7,1.985){\line(1,0){.3}}
\put(2.7,2.015){\line(1,0){.3}}
\put(2,1.985){\line(1,0){.3}}
\put(3,2){\line(1,0){1}}
\put(2,2){\circle*{.06}}
\put(2.01,1.98){\line(1,2){.15}}
\put(1.98,2){\line(1,2){.15}}
\put(2.185,2.33){\line(1,2){.15}}
\put(2.155,2.35){\line(1,2){.15}}
\put(2.36,2.68){\line(1,2){.15}}
\put(2.33,2.7){\line(1,2){.15}}
\put(2,2){\line(1,-2){.5}}
\put(3,2){\circle*{.06}}
\put(2.845,2.35){\line(-1,2){.15}}
\put(2.815,2.33){\line(-1,2){.15}}
\put(2.67,2.7){\line(-1,2){.15}}
\put(2.64,2.68){\line(-1,2){.15}}
\put(3.02,2){\line(-1,2){.15}}
\put(2.99,1.98){\line(-1,2){.15}}
\put(3,2){\line(-1,-2){.5}}
\put(4,2){\circle*{.06}}
\put(2.5,3){\circle*{.06}}
\put(2.5,3){\line(0,-1){.95}}
\put(2.5,1){\circle*{.06}}
\put(2.5,1){\line(0,1){.95}}
\end{picture}\\[-16mm]
\end{cor}

\proof
This follows directly from the fact that the subset $\D'$ has six elements of
this kind, since it stems from the complete vanishing lattice $L',\D'$.
\qed

We get a extended result for $K3$ surfaces analogously.

\begin{lemma}
\labell{K3-pb}
Let $X$ be an elliptic $K3$ surface then there is a set $\D$ of classes of $L_X$
such that the elements of $\D$ generate $L_X$ and are conjugated under the
action of $\G_\D$ and such that $\D$ contains six elements with intersection
pattern as in the corollary above.
\end{lemma}

\proof
Let $\D$ be the union of classes as given in prop.\,\ref{v-pb} with the pull back
of the classes $\a,\a',\s$ of an elliptic $K3$ surface as in prop.\,\ref{K3ex}
along a deformation equivalence through elliptic surfaces, then the first and
third property are obvious from the above corollaries, and for the second we argue
as in cor.\ref{D-conj} with the additional classes.
\qed

\begin{thm}
\label{final}
Let $X$ be a minimal regular elliptic surface of positive geometric genus and
canonical class $k_X$, then the group $\G_X$ of orthogonal transformations of its
homology lattice
$L_X$ generated by the monodromy groups of all families containing $X$ is
$$
O'_k(L_X):=\{\cg\in O(L_X)|\cg(k_X)=k_X,\cg \text{ has positive real spinor
norm}\}
$$
\end{thm}

\proof
Since any monodromy element fixes the canonical class and is induced from an
orientation preserving diffeomorphism we infer from \cite{kml} that there is an
inclusion $\G_X\inj\Ospk(L_X)$ of subgroups of the orthogonal group of $L_X$.

Define $\D_X$ to be the set of all classes of $L_X$ orthogonal to the fibre class
which are represented by vanishing cycles in a family containing $X$ such that the
associated reflection $\cg$ is in the monodromy group. Hence
$\G_{\D_X,L_X}$ is a subgroup of $\G_X$.

We denote the orthogonal complement of $k_X$ in $L_X$ by $L_k$, and observe that
$\D_X$ is a subset of $L_k$.
Thus $L_k$ is stabilized by $\G_{\D_X,L_X}$ and we get a surjection
$\G_{\D_X,L_X}\to\!\!\!\!\!\to\G_{\D_X,L_k}$ with the latter a subgroup of
$O'(L_k)$.

Since elements of $\Ospk(L_X)$ fix not only the canonical class but also its
orthogonal complement $L_k$ and since the spinor norm is preserved under the
restriction to the action on $L_k$ we get a map $r_k:\Ospk(L_X)\to O'(L_k)$.

This map $r_k$ is the identity map in case of a $K3$ surface. Otherwise
$k_X$ is a non-trivial multiple of the fibre class and $L_k$ coincides with $L$.
With $\cg_0\in\Ospk$ a preimage of the identity in $O'(L)$ we get
$$
\begin{array}{llcl}
\cg_0(c).c'=c.c' & \forall c\in L_X,c'\in L & \Rightarrow & \cg_0(c)-c\in L^\perp\\
\cg_0(c).k_X=c.k_X & \forall c\in L_X & \Rightarrow & \cg_0(c)-c\in L\\
\multicolumn{2}{l}{\cg_0(c)^2=c^2+2c.(\cg_0(c)-c)+(\cg_0(c)-c)^2}.
\end{array}
$$
First we conclude that $\cg_0(c)-c$ is contained in the radical of $L$, then we
infer that $\cg_0(c)=c$ since either $c$ is in $L$ or $c$ defines an injective
form on the radical. Thus $r_L$ is injective in either case and we get a
commutative diagram
$$
\begin{array}{ccccc}
\G_X & \subset & \Ospk(L_X) & \inj & O'(L_k)\\
\cup &&&&\parallel\\
\G_{\D_X,L_X} & \to\!\!\!\!\!\to & \G_{\D_X,L_k} & \subset & O'(L_k).
\end{array}
$$
To apply the criterion of Ebeling \cite{habil} in the case that $X$ is not a $K3$
surface, let $\D$ be as in prop.\,\ref{v-pb} and observe that
$$
\D_k:=\G_\D\cdot\D
$$
is a subset of $L_k$ of classes of square $-2$ such that $L_k,\D_k$ is a complete
vanishing lattice.
This is immediate from the corollaries to the proposition since $\G_\D=\G_{\D_k}$
due to the fact that
$$
\cg_{\cg(\a)}=\cg\circ\cg_\a\circ\cg\inv
$$
for $\cg\in\G_\D$ and $\cg_{\cg(\a)},\cg_\a$ the reflections associated to
$\cg(\a)\in\D_k,\a\in\D$.
We may conclude that $\G_{\D_k}=O'(L_k)$ and that $\D\subset\D_X\subset\D_k$.\\
In the case that $X$ is a $K3$ surface we take $\D$ from lemma \ref{K3-pb} and
notice that for
$$
\D_k:=\G_\D\cdot\D
$$
we may conclude as above that $\G_{\D_k}=O'(L_k)$ and that
$\D\subset\D_X\subset\D_k$.\\
Thus in either case even $\G_{\D_X,L}=O'(L)$, hence the maps in the top
row of the commutative diagram must be isomorphisms and $\G_X=\Ospk(L_X)$ is proved.
\qed

Of course the conclusion holds for arbitrary $K3$ surfaces by the deformation
equivalence of all $K3$ surface and has been proved before as a special case of a
complete intersection surface with positive geometric genus, cf.\ \cite{habil}.\\
To compare the monodromy action with the action of the diffeomorphism group we
have to invoke results of \cite{kml} to get

\begin{cor}
Let $X$ be as above then $\G_X$ is the stabilizer group of the canonical class in
the image of the natural representation of the group of orientation preserving
diffeomorphisms on $L_X$ and of index two.
\end{cor}

\end{document}